\documentclass[10pt]{amsart}
\usepackage{graphicx}
\usepackage{hyperref}
\usepackage{amssymb}
\usepackage[dvipsnames]{xcolor}
\usepackage{tkz-graph}
\usepackage[utf8]{inputenc}

\usepackage{soulutf8}
\usepackage{indentfirst}

\makeatletter
\@namedef{subjclassname@2010}{%
\textup{2010} Mathematics Subject Classification}
\makeatother
\newtheorem{thm}{Theorem}[section]
\newtheorem{cor}[thm]{Corollary}
\newtheorem{lem}[thm]{Lemma}
\newtheorem{prop}[thm]{Proposition}

\theoremstyle{definition}

\newtheorem{rem}[thm]{Remark}
\newtheorem{exa}[thm]{Example}
\newtheorem{exe}[thm]{Exercise}
\newtheorem{ass}[thm]{Assumption}

\newcommand{\numberset}{\mathbb}
\newcommand{\R}{\numberset{R}}

\GraphInit[vstyle = Shade]
\tikzset{
  LabelStyle/.style = { rectangle, rounded corners, draw,
                        minimum width = 2em, fill = Green,%yellow!50,
                        text = red, font = \bfseries },
  VertexStyle/.append style = { inner sep=5pt,
                                font = \normalsize\bfseries},
  EdgeStyle/.append style = {->, bend left} }

\begin{document}

\baselineskip=17pt
\title[Sample paper]{Geometric Proofs and algebraic functions}
\author[D. Maran]{Davide Antonio Nello Maran}
%\email{davide.maran@mail.polimi.it}
\date{\today}
\maketitle
\begin{abstract}
In this paper I present a kind of proof for classical Euclidean geometric problems which relies on both synthetic and analytic geometry. Using the elementary tools of polynomial algebra and multivariate calculus we manage to reduce the proof of some general theorems to a finite number of particular cases.
\end{abstract}

\section{Before starting}
For 20 centuries "geometry" had been a synonymous for "Euclidean Geometry", a refined form of math where the proof of the theorems where all oriented in apply the famous \emph{5 postulates of Euclid} in endless different ways. Demonstrations of this type, now called "synthetic", relied mostly on the creativity of the author, who should be able to find the key of every problem without fixed strategies.

At the middle of the $17^{th}$ century the world saw the dawn of new method of interpret geometry: the work of René Descartes and Pierre de Fermat led to the possibility of representing every point with two or more scalar coordinates. In this vision, lines, circles and every other object becomes an equation in the two coordinates which is satisfied by the points which belong to it. Following this pact, every proposition becomes an equation in some real variables which correspond to the cases where the proposition is true, so that every theorem becomes an identity.

To use of this powerful instrument, we pay a price in terms of beauty of the demonstrations: with the analytic approach every problem of proving a proposition brings to an oceanic amount of computations to establish if the equation associated to that proposition is always satified. 

What I am doing in this paper is to find an hybrid approach that tries to climb these two pillars of human knowing at the same time, providing proofs which are not as elegant as the ones of synthetic geometry neither as plain as the one of analytic one, but show the simplest sides of both approaches, using the synthetic one for "test cases" and the analytic one to say how many particular cases I need to have a proof which is always valid. In fact, if we see every proof in the analytic geometry as the verification of
$$f(a,b,c,\dots)=0$$
for any chose of the variables, my approach, which i call \textbf{Autarkic} achieves the result following an unusual path which is explained in the following scheme

\begin{align} \text{part. cases} & \xleftarrow[simpler\ proof]{} \text{general case}\\ \updownarrow & \qquad \qquad \qquad \updownarrow{\textcolor{Red}{\text{analytic\ geometry}}}  \\ P_{i=1,...}(a,b,c,..)=0 & \xrightarrow[use\ thm.\ \ \ref{main}]{} f(a,b,c,\dots)=0 \end{align}

starting from the general case of a proposition i.e. "\textit{for every triangle...}", instead of directly finding its associated function $f(a,b,c...)$, which is the method of cartesian geometry and would be very onerous, we restrict our proof to some particular cases, for example "\emph{for every rectangular triangle...}" where you could easy end the proof in an Euclidean way. At this point we will see how to associate some algebraic functions $P_i$ to these cases and to apply a non-trivial result that allows us to deduce the fact that $f=0$ everywhere.

\subsection{The main theorem}
It is unusual that the proof of the main theorem of an article comes right after the introduction. Well, the following is the result on which I am basing all my work, but at the same time is not the core of this work, since its importance on geometry will be seen much later in the paper.
\begin{thm}
Let $f: \R^n\to \R$ be a polynomial function of degree $g$ and $P: \R^n\to \R$ another $n-$variate polynomial such that
$$\forall \vec a\in \R^n\ P(\vec a)\implies f(\vec a)=0$$
if there exist a line $\vec r: \R \to \R^n$ such that $P(\vec r(k))$ has $g$ disjoint single roots, then $f(\vec a)=0$ everywhere.
\label{main}
\end{thm}

\begin{proof}
Indeed we can represent the line in $\R^n$ as a function $\R \to \R^n$ of this form
$$\vec r(k;\vec b_0)=\vec a_0k+\vec b_0$$
Now let
$$p(k,\vec b)=P(\vec r(k;\vec b))=P(\vec a_0k+\vec b)$$
this function goes from $\R \times \R^n$ to $\R$ and we know that
\begin{itemize}
\item $p(\cdot, \vec b)$ is a polynomial function for all $\vec b$.
\item $p(k,\cdot)$ is a polynomial function for every $k$
\item $p(k,\vec b_0)$ has $g$ distinct roots $k_1>k_2>k_3...>k_g$
\end{itemize}
moreover, if for some $i$ $$\frac{\partial p}{\partial k}(k_i,\vec b_0)=0$$
we would have that $k_0$ is a double root of $p(k, \vec b_0)$\footnote{it is a well-known fact that if a root of a polynomial is also a maximum (or a minimum) then it is a multiple root}
so we can say that under our assumptions
$$\forall i=1,...g\ \qquad \frac{\partial p}{\partial k}(k_i,\vec b_0)\neq 0$$
under this fact we can apply the implicit function theorem so that in a neighbourhood $B$ of $\vec b_0$ there are $g$ functions $s_1,s_2,...s_g:$ $\R^n\to \R$ such that
$$\forall \vec b\in B,\ \ \  s_1(\vec b)>s_1(\vec b)>...>s_n(\vec b)$$
$$\forall i=1,...g\ \qquad \ p(s_i(\vec b),\vec b)=0$$
Let us consider now, for every $\vec b\in B$ fixed,
$$f(r(k;\vec b))$$
we know that this polynomial function of $k$ has degree $g-1$ and takes value zero in at least $g$ points, which we have just called
$$s_1(\vec b),s_2(\vec b),...s_g(\vec b)$$
this is possible only if $f(r(k;\vec b))$ is zero for every $k$!

So, we have in particular that, evaluating the function in $k=0$,
$$0=f(r(0;\vec b))=f(\vec b)$$
that holds for every $\vec b\in B$. But $B\subset \R^n$ is an open set by construction, so $f$ has value zero on a whole open set. Since $f$ is a polynomial function it entails that $f=0$ everywhere
\end{proof}

%% 	RIMARCARE. ORA!

\begin{rem}
The assumption of single roots is fundamental. Indeed take $P(a,b)=(a+b)^3$ and $f(a,b)=a^2-b^2$. Of course if
$$0=(a+b)^3=P(a,b)$$
then $a=-b$, so
$$f(a,b)=a^2-b^2=0$$
and taking $\vec r: \R \to \R^2$ as
$$\vec r(k)=\begin{pmatrix} 1\\ k\end{pmatrix}$$
the equation $f(\vec r(k))=0$ becomes
\begin{equation}k^3+3k^2+3k+1=0 \label{terzo}\end{equation}
which has three solution (fundamental theorem of algebra). Nevertheless, $f$ is \emph{not} identically zero!\\
This happens since the roots of \eqref{terzo} are not distinct, in fact it is straightforward to verify that all the roots correspond to $-1$.
\end{rem}

%%		corollarino
\section{A first way to apply the main theorem}
As first example of what can we do with our main theorem, we present a naive, coordinate-free approach based on represent a triangle in terms of the lengths of its sides and every property as a function of them.

\subsection{When the notable angles are enough?}

In order to apply our main theorem we have to find some polynomial equation which is relevant to geometry. There can be made endless possible choices and in this section we are going to focus on one of the most famous. 
Let $a,b$ and $c$ be the lengths of the sides $BC,AC$ and $AB$ of a triangle and recall that, from the so called "law of cosines", we have, calling $\gamma$ the angle at the opposite to $AB$,
\begin{equation}a^2+b^2-2\cos(\gamma) ab=c^2 \label{cos}\end{equation}
This means that, fixing $\gamma$, the polynomial 
$$P(a,b,c)=a^2+b^2-2\cos(\gamma) ab-c^2$$
is zero if and only if the triangle $ABC$ has the angle $\widehat {ACB}=\gamma$. It is important to remember this interpretation while proving this result:

\begin{cor}
%%%%%%%%%%%%%%%%%%%%%
 \label{coroll}
%%%%%%%%%%%%%%%%%%%%%
Let $f:\R^3 \to \R$ be a polynomial function of degree $2g-1$. If $f$ gives value zero when evaluated on the sets
$$S_i:=\{a,b,c: a^2+b^2-c^2=2\lambda_i ab\}$$
for at least $g$ distinct values of $\lambda_i\in (-1,1)$ then $f$ is identically zero.
\end{cor}
\begin{proof}
In order to use the previous theorem, first of all we need a line. Let 
$$\vec r(k)=\begin{pmatrix}1\\1\\1+k\end{pmatrix}$$
and of course, we need a polynomial, say
$$P(a,b,c)=\prod_{i=1}^g (a^2+b^2-c^2-2\lambda_i ab)=$$
$$=(a^2+b^2-c^2-2\lambda_1 ab)(a^2+b^2-c^2-2\lambda_2 ab)...(a^2+b^2-c^2-2\lambda_g ab)$$
Now, which are the roots of
$$P(\vec r(k))?$$
of course, it is enough that one term of the product is zero to make the whole product be zero!
Thus, for any $i=1...g$ the numbers
$$k_{i,1}=\sqrt{2-2\lambda_i}-1$$
$$k_{i,2}=-\sqrt{2-2\lambda_i}-1$$
are roots of $P(\vec r(k))$, and it is straightforward to verify that distinct values of $i$ correspond distinct values of $k_{i,1},k_{i,2}$ (I could say that it is a useful exercise for the reader but I would prefer you not to burn this paper).

So, we have that $P(\vec r(k))$ has $2g$ distinct roots and all of them have multeplicity one (the polynomial has degree $2g$ so if the roots were multiple, they where less than $2g$). Since by assumption
$$P(a,b,c)=0\implies f(a,b,c)=0$$
we can apply \ref{main} to end the proof.
\end{proof}

%%	 dis	corsivo

However, finding an interpretation for the sets $S_i$ of the previous corollary \ref{coroll} is useless until we find an interpretation also for the polynomial function $f$. What I am going to convince you, which is the main assumption under this work is that

\begin{ass}\label{ass}
For every property of a triangle in the Euclidean geometry, there is polynomial a function $f: \R^3\to \R$ such that $f(a,b,c)=0$ if and only if the property is valid for a triangle $ABC$ with $AB=c,\ BC=a,\ CA=b$
\end{ass}

For now, we cannot give a rigorous value to this assumption, since we have not well defined what we mean for "property"! However, in the next section \ref{cartesio} about the cartesian geometry this assumption will become clearer and will be supported by one of the most famous result in Galois-algebra. 
For now, I just want you to consider that thanks to the third congruence criterion the sides of a triangle determine in an unique way all the geometric features of the triangle itself so that anything can be constructed starting from a triangle is indeed a function of the three sides.
Often this function is not a polynomial, but as we are going to see, there are several cases in which this quantity is defined by a polynomial of $a,b,c$ in implicit form\footnote{i.e. $f(a,b,c)$ is such that exists polynomial $P(a,b,c,f(a,b,c))=0$ for any $a,b,c$}.

%% discorso

Let us see these examples:
\begin{itemize}
\item Thanks to Erone's formula the area of the triangle satisfies
$$ 4Area^2={(a+b+c)(a+b-c)(a-b+c)(-a+b+c)}$$
which is an equation of fourth degree in $a,b,c$ that allows to express the aera of the triangle as an implicit polynomial function of the sides.
\item Calling $h_a,h_b,h_c$ the length of the heights relative to $A,B,C$ respectively, we know from previous point that they satisfy
$$ (h_aa)^2={(a+b+c)(a+b-c)(a-b+c)(-a+b+c)}$$
$$ (h_bb)^2={(a+b+c)(a+b-c)(a-b+c)(-a+b+c)}$$
$$ (h_cc)^2={(a+b+c)(a+b-c)(a-b+c)(-a+b+c)}$$
\item The distance $d$ of the incenter from the sides can be computed as well
$$r^2(a+b+c)^2={(a+b+c)(a+b-c)(a-b+c)(-a+b+c)}$$
which can be simplified into
$$r^2(a+b+c)={(a+b-c)(a-b+c)(-a+b+c)}$$

\end{itemize}

Using this consideration, for every geometrical proposition about some of these objects we can find at least one polynomial function into the variables $a,b,c$ that takes value zero if and only if a triangle $ABC$ of sides $BC=a,CA=b,AB=c$ satisfy the proposition.

\begin{exa}Consider the proposition\\
\label{AA}
\emph{(P):\ \ The height relative to the side $AB$ has 3 times the length of the distance from the incenter}

\begin{figure}[ht]
  \includegraphics[scale=0.6]{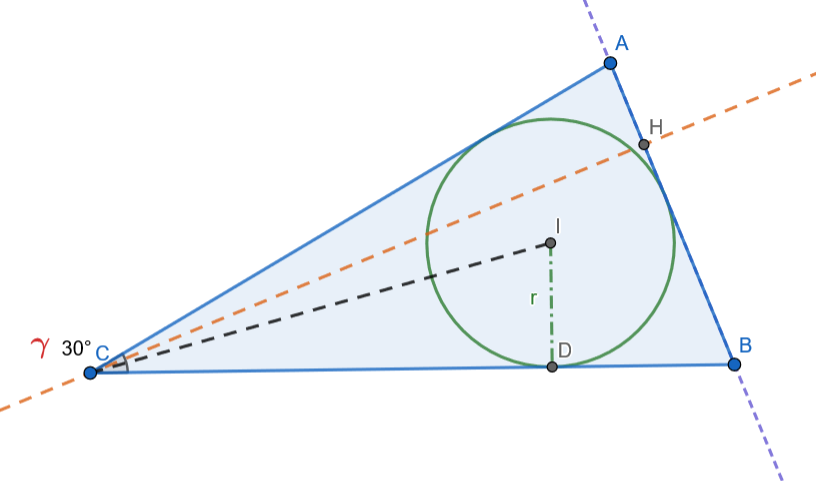}
  \caption{triangle of example \ref{AA}}
  \label{A}
\end{figure}
we have just seen that
$$ (h_aa)^2={(a+b+c)(a+b-c)(a-b+c)(-a+b+c)}$$
$$r^2(a+b+c)={(a+b-c)(a-b+c)(-a+b+c)}$$
so we know that (P) holds if and only if
$$\sqrt{\frac{3(a+b-c)(a-b+c)(-a+b+c)}{(a+b+c)}}=\sqrt{\frac{(a+b+c)(a+b-c)(a-b+c)(-a+b+c)}{a^2}}$$
which is not a polynomial function since it contains the square roots and the denominators but, squaring both members,

$$\frac{9(a+b-c)(a-b+c)(-a+b+c)}{(a+b+c)}=\frac{(a+b+c)(a+b-c)(a-b+c)(-a+b+c)}{a^2}$$
and reordering the terms, we can define a polynomial function 
$$f(a,b,c):=9(a+b-c)(a-b+c)(-a+b+c)a^2-\dots \qquad \qquad \qquad$$
\begin{equation}\qquad \qquad \qquad \dots+(a+b+c)(a+b-c)(a-b+c)(-a+b+c)(a+b+c)\label{brutto} \end{equation}
which is zero if and only if a triangle of sides $(a,b,c)$ satisfies $(P)$.
\end{exa}

The property (P) is satisfied if the lengths of the sides of the triangle are a root of the polynomial function \eqref{brutto}, so how can we know if it is valid for all the triangles or not?

\emph{Easy! \ (P) holds for every triangle if and only if \eqref{brutto} is identically zero!}

The problem is that \eqref{brutto} is a 5 degree polynomial nearly as long as a \LaTeX $ \ $ document can support but much more than what human brain can, so we want to fing an easier way to know whether (P) is a theorem or not. 

For now we can tell for sure that $f(a,b,c)$ is a polynomial functions of degree $5$, and since $5=2\times 3-1$ we know for \ref{coroll} that $f$ is always zero if it is zero on at least three sets of the form
$$S_i:=\{a,b,c:\quad a^2+b^2-c^2=2\lambda_i ab\}$$
for different values of $\lambda_i$. Note that, as stated by \ref{cos}, by the analogy with the cosines law, to prove that $f$ is zero on these sets means exactly to prove that the theorem holds for the triangles $ABC$ such that the angle $\gamma$ between sides $CA$ (of length $b$) and $BC$ (of length $a$) is fixed with a value $\gamma_i=\arccos{(\lambda_i)}$.

Note that we can choose the three values of $\lambda_i$ (and so $\gamma_i$) however we like! For example, we could use $\lambda_i\in \{\frac{\sqrt 3}{2},\frac{1}{2},0\}$. In these cases, we find ourself with particular triangles with an angle of $\pi/6,\pi/3$ or $\pi/2$, so that we can do the verification by hand.

Let us start with the triangle having $\gamma=\pi/6$.

In this case, it is easy to say that the proposition may not hold: suppose that the triangle is convex, then the distance between the incenter $I$ and $C$ is more than $\overline{CH}$ (the length of the relative height). $\overline{CH}$ is linked to the distance $r$ of the incenter to the sides \footnote{remember that by definition the distance of the incenter to the three sides is equal} by a simple relation
$$r=\overline{CI}\sin \biggl(\frac{\pi/6}{2}\biggr)$$
since $\sin (\pi/6)<1/3$ we have
$$\overline{CH}>\overline{CI}>3r$$
that shows that (P) may fail for triangles with $\gamma=\pi/6$. It is not even necessary to consider the other two cases to say that (P) is not a theorem.

\subsection{Have we really used our new result?}
In the previous example, corollary \ref{coroll} is not actually used! The fact that if a statement is false we can disprove it by only showing a counter example is a basic logical fact that is pretty well know. In fact we were trying to use \ref{coroll} when we found that even for the first angle the thesis was not true!
Actually, the right path to face any geometrical proof with corollary \ref{coroll} is
\begin{enumerate}\label{procedura}
\item Find the degree $g$ of the polynomial function associated to the problem (actually this is often the most difficult part)
\item Chose $g/2+1$ notable angles such that the statement becomes much simpler for a triangle with one of those angles
\item Try to prove the theorem for these $g/2+1$ special cases. If at a certain point we arrive to a contradiction then the general theorem is false. If we manage to end this procedure the theorem is true.
\end{enumerate} 
%So ou goal is to provide an example where we can prove a real theorem with our method. The more we can even avoid using \ref{coroll} and try to open our mind to use \ref{main} in an even more efficient way.
this approach, that I call \emph{Autarkic}, can be summarized in the following scheme:
If we indicate the polynomial function $f$ associated to the problem and $P_1,P_2,\dots$ the polynomials of the form 
$$P_i=a^2+b^2-c^2=2\lambda_i ab$$
that represent the special cases of a triangle with a fixed angle we proceed like that, starting from (1).

\begin{align} \text{part. cases (2)} & \xleftarrow[simpler\ proof]{} \text{general case (1)}\\ \updownarrow & \qquad \qquad \qquad \updownarrow  \\ P_1,P_2,\dots (3)& \xrightarrow[use\ cor.\ \ \ref{coroll}]{} f(a,b,c) \end{align}

Instead of computing $f$ explicitly and show whether it is zero or not, we pass trough the special cases of some fixed angles.
%%lo rifaccio sotto o o sopra? sopra!!!

\section{The cartesian coordinate version}\label{cartesio}
After seeing the previous naive method, we have two goals: find an example of famous result that can be proved with our algebraic-oriented autakic method \ref{main} and to give a rigorous meaning to the assumption \ref{ass}.

For the first issue let us try to find an autarkic proof of one wonderful theorem that has made the history of geometry.
\begin{thm}\label{Euler}
\textbf{(Euler's line)} For every every triangle $ABC$ there exists a line which passes through the centroid, the orthocenter and the circumcenter
\end{thm}
This wonderful theorem has is still a source of inspiration for research. For example, just four years ago, \cite{Hun} showed a property of this line in a triangle built on the golden rectangle!

In order to prove this theorem, we would like to find, as in the previous section, a polynomial function in the lengths $a,b,c$ of the sides of a triangle that has value zero if and only if the three notable points are aligned. As we have seen, we can deduce the heights of a triangle from its sides and it can be shown that we could also write a polynomial in implicit form that links the position of the three notable points to $(a,b,c)$. However, this procedure is very difficult and our strategy requires that the degree of the polynomial that comes out is as low as possible, in order to be able to find a reasonable number of special cases.

For this reason, we do a re-parametrization: we will write the polynomial function associated to the problem in three variables which represent the triangle but that are not the sides. To do this, we make use of one of the most succesful mathematical tools of all time: \emph{the cartesian coordinate system.}

Using that, any triangle can be parametrized with three variables in many ways.
\begin{figure}[ht]
  \includegraphics[scale=0.8]{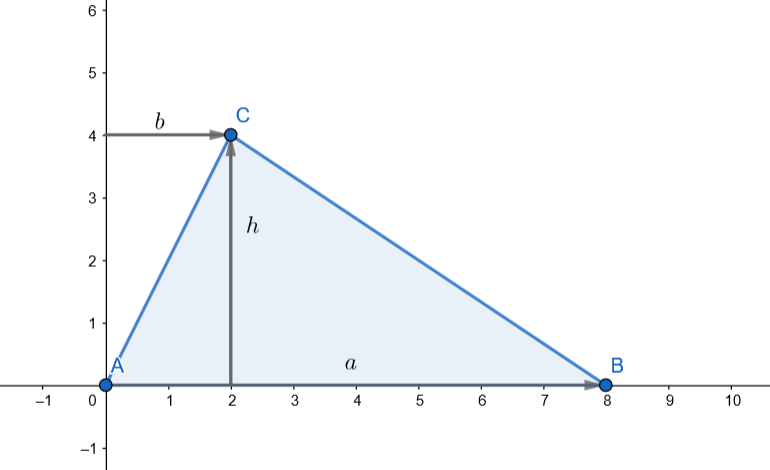}
  \caption{our new variables $a,b,h$ which describe a triangle in the coordinate system}
  \label{A}
\end{figure}
For example, we will use the three variables $a:$ x-coordinate of $B$, $b:$ x-coordinate of $C$ and $h:$ y-coordinate of $C$, imposing that the side $AB$ lies on the $x-axis$. Summarizing we have
\begin{itemize}
\item $A=(0,0)$
\item $B=(a,0)$
\item $C=(b,h)$
\end{itemize}

\begin{rem}
Of course, our new three variables $a,b,c$ are linked to the previous ones $a',b',c'$ which represented the length of the sides. Indeed,
$$a=c',\qquad h^2+b^2=b'^2$$
$$(a-b)^2+h^2=a'^2$$
so that %$$(a-a')(h^2+b^2-b'^2)((a-b)^2+h^2-c'^2)=0$$
everything that can be expressed as implicit polynomial function of $a',b',c'$ can be also expressed in terms of $a,b,h$.
\end{rem}

To use this coordinate system allows us to finally prove our \ref{ass}!

Defining a \emph{constructive step} as an elemental construction that is allowed under the axioms of Euclidean geometry, it has been proved that
\begin{thm}
Given some points $A_{1,...n}$ in the cartesian coordinate system of coordinates $(a_{i,x},a_{i,y})$, every \emph{constructive step} leads to a point $A_{n+1}=(a_{n+1,x},a_{n+1,y})$ such that there exist polynomial functions $P_x,P_y,Q_x,Q_y$
$$a_{n+1,x}=P_x(a_{i,x},a_{i,y})+\sqrt{Q_x(a_{i,x},a_{i,y})}$$
$$a_{n+1,y}=P_y(a_{i,x},a_{i,y})+\sqrt{Q_y(a_{i,x},a_{i,y})}$$

\end{thm}
Iterating this theorem, we find an important corollary
\begin{cor}
All the points which can be constructed from a given set of points in a finite number of steps have coordinates which are an implicit polynomial function of the initial coordinates
\end{cor}
The converse is not true, starting from some points of integer coordinates we won't ever be able to build points like $(0,\sqrt[3]{2})$ although the coordinates  solve a polynomial equation of integer coefficients\footnote{to see a real application with GeoGebra see \cite{Der}}. 
The theoretic foundations of this result are very complicated and to understand them, I suggest \cite{Hol}.
This means that every point that can be constructed in Euclidean geometry has coordinates that solve a polynomial equation with the original coordinates as parameters\footnote{\href{https://math.stackexchange.com/questions/3612233/a-polynomial-relation-between-the-lengths-of-segments-that-are-constructible-wit}{there are also discussions open on the topic}}.

This last corollary encourages us to focus our attention on the problem parametrized with the cartesian coordinate system. What we are going to do is to change 
a bit our procedure \ref{procedura}. We could see that while dealing with analytic geometry, to translate the relation between two objects into an equation is much simpler (that is the reason why Cartesio invented it), and indeed it is much simpler even estimating the degrees of the equation that comes out. Unfortunately, the results which make these estimations are very boring and unnecessarily complicated, so I preferred to leave them in the appendix \ref{boring}.

After this change of environment, instead of using as special cases the triangles with one angle of $\pi/6, \pi/4, \pi/2$ and apply \ref{coroll} we are going to try with other classes of "special triangles". The procedure to solve the problem $(P)$ of stating if a given property is satisfied from all the triangles changes as follows

\begin{enumerate}
\item Find the degree $g$ of the polynomial function associated to the problem \textcolor{Green}{with the result found in Appendix \ref{boring}}
\item Choose \st{$g/2+1$ notable angles such that the statement becomes much simpler for a triangle with one of those angles} \textcolor{Green}{?? "special triangles" which make the proof easier}
\item Try to prove the theorem for these \st{$g/2+1$} \textcolor{Green}{??} special cases. If at a certain point we arrive to a contradiction then the general theorem is false. If we manage to end this procedure the theorem is true.
\end{enumerate} 

I have voluntarily left the question dots instead of telling you the triangles that we are using because, despite I have anticipated that we will test our strategy on the Euler's line problem, it could be used for any geometrical problem and in every situation, so it is important for every problem to be able to find the special cases which make the proof easier. By the way, we are ready to face the Euler's line.

\subsection{Example of analytic use: Euler's line}

This is how the problem looks for any triangle with the parametrization previously defined.

\begin{figure}[ht]
  \includegraphics[scale=0.8]{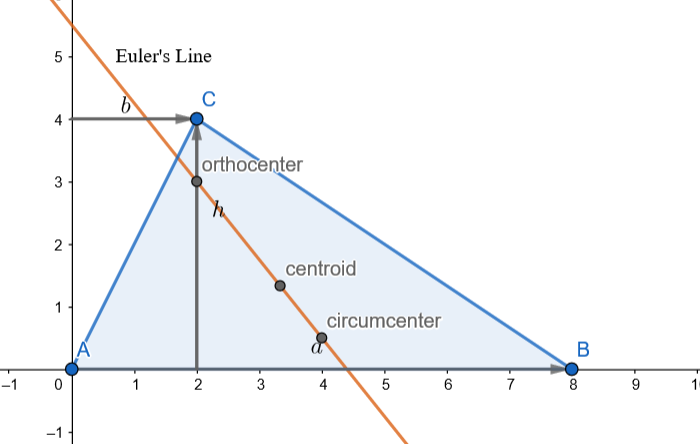}
  %\caption{Per la parte informatica ho chiesto l'aiuto di un hacker di prim'ordine}
  \label{P}
\end{figure}

Now, the coordinates of the three notable points are rational functions of $a,b,h$, we have to estimate their degree. 

\begin{itemize}
\item The centroid is the simplest, its coordinates are just the mean of the coordinates of the three vertexes. The computation is so simple we do not actually need to make any estimation and provide $$G=\biggl(\frac{a+b}{3},\frac{h}{3} \biggr)$$
so that both coordinated are polynomial functions of degree $1$\\

\item The otrhocenter requires some extra computations, so in order to avoid them it is enough to use the lemmas of appendix \ref{boring}.
The orthocenter is the intersection of the height passing from $C$ and the height passing from $B$. In particular, the height $\vec q$ relative to $C$ is simply a vertical axis with $x=b$, while 
%so that it can be written as
%$$\vec q(j)=\begin{pmatrix} b\\ j\end{pmatrix}$$
moving to the height $\vec r$ passing from B we know that, being orthogonal to $AC$, its equation must be
$$\vec r(k)=\begin{pmatrix} -h \\ b \end{pmatrix} k+\begin{pmatrix} a\\ 0\end{pmatrix}$$
Since the components of $\vec q,\ \vec r$ are linear in $a,b,h$, we can apply \ref{2Rpt} to prove that the $x-$coordinate of their intersection point, the orthocenter, has at most degree $2$ over $1$ 
\\

\item For the circumcenter the situation is not very different. It is the point of intersection of the axes of $AB$ and of $AC$, so we have the coefficient and one point passing in both cases.

Axis $\vec p$ of $AB$ is vertical and passes for the midpoint of $AB$, so
$$\vec p(j)=\begin{pmatrix} a/2 \\ j\end{pmatrix}$$
while the axis $\vec s$ of $AC$ is orthogonal to $AC$ and passes from the midpoint of $AC$, so its equation is
$$\vec s(k)=\begin{pmatrix} -h \\ b \end{pmatrix} k+\begin{pmatrix} b/2\\ h/2\end{pmatrix}$$
again, since the components of $\vec p,\vec s$ are linear in $a,b,c$, lemma \ref{2Rpt} ensures that the $y-$coordinate of their intersection point, the circumcenter, has at most degree $2$ on $1$ 
\end{itemize}
With these considerations, we have estimated the degrees of the coordinates of the centroid $G$, the hortocenter $H$ and the circumcenter $Q$ as reported in the following table. 
\begin{table}[ht]
%\caption{supporto e calco}
\center
\label{tab:buoni}
\begin{tabular}{|c|c|c|}
\hline
Point & x-coordinate & y-coordinate \\
\hline
G & $ 1$ & $1$ \\
H & $ 1$ & $2/1$\\
Q & $ 1$ & $2/1$\\

\hline
\end{tabular}
\end{table}

Now we have to translate the property \emph{"G,H,Q are collinear"} into a real function $f(G(a,b,h),H(a,b,h),Q(a,b,h))$, so that substituting the functions which express the coordinates of $G,H,Q$ into that we get $0$ if and only if the triangle corresponding ro $(a,b,h)$ satisfies this property. 

Knowing that any three points $(d_x,d_y),\ (e_x,e_y)\ (h_x,h_y)$ are collinear if and only if they satisfy
\begin{equation}\label{coll}(e_x-d_x)(h_y-d_y)=(h_x-d_x)(e_y-d_y)\end{equation}

we can define $f(a,b,h)$ by substituting the coordinates of $G,H,Q$ respectively. However, doing it explicitly would be a waste of time! with our method  is enough to find the degree of $f$ from the degrees of the coordinates of $G,H,Q$:
 
$$ (1^{st} -1^{st})\biggl ( \frac{2^{nd}}{1^{st}}-1^{st} \biggr) =(1^{st} -1^{st})\biggl ( \frac{2^{nd}}{1^{st}}-1^{st} \biggr)$$
which is, after doing the multiplication,
$$ \biggl ( \frac{3^{rd}-3^{rd}}{1^{st}} \biggr) = \biggl ( \frac{3^{rd}-3^{rd}}{1^{st}}\biggr)$$
that becomes
$$4^{th}=4^{th}$$
In the end we have estimated that there exists an equation of degree $4$ into $a,b,h$ such that it is satisfied if and only if the triangle with parameters given by $a_0,b_0,h_0$ admits the Euler line, so we can define $f(a,b,c)$ simply as the difference between the left-hand side and the right-hand side; its maximum degree is still four.

It has been a long journey to arrive to that result, but only because I had to do all the steps to explain the procedure! In fact, it is easier to make it that to tell it!

Anyway, we now seek to prove that the Euler line exists for every triangle, and in order to do that we have to prove that $f$ is identically zero by exploiting the fact that it has almost degree $4$ and it is zero for some particular values of $(a,b,h)$ which correspond exactly to the classes of triangles that we are going to see in the following and last step.

\subsection{Last step: particular cases}
We have finally arrived to the only part of real demonstration, as intended in the Euclidean sense, and finally we have to choose some triangles where the existence of Euler's line is trivial.
Here it is useless to consider special cases of triangles with one angle fixed, since the proof in that cases becomes not simpler than what it is in the general case, but we can indeed manipulate the lengths of the sides.
For example, if the triangle $ABC$ is \emph{isosceles} the existence of Euler line is trivial: all the three points lay on the axis by symmetry!

This mean that whenever $AB=BC$, $BC=CA$ or $CA=AB$ we have our result. But what does it imply in terms of the three variables $a,b,h$ we have used to parametrize the triangle? 

Using Pithagoras' theorem we have
$$AB=BC \iff a^2=b^2+h^2$$
\begin{equation} BC=CA \iff a=2b \label{iso}\end{equation}
$$CA=AB \iff h^2+(a-b)^2=a^2 \iff h^2+b^2=2ab$$

\begin{figure}[ht]
  \includegraphics[scale=0.8]{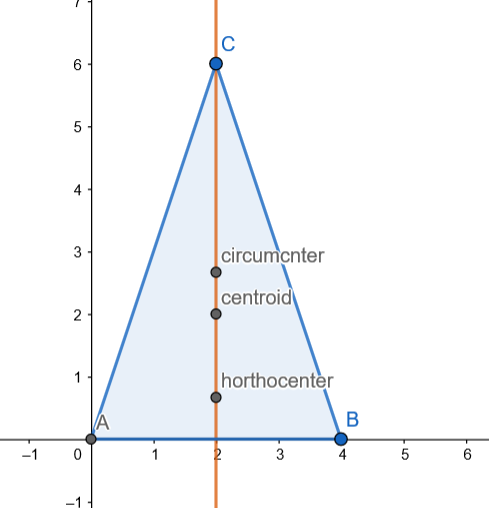}
  \caption{Euler's line coincides with the axis of symmetry for an isosceles triangle}
  \label{P}
\end{figure}

so that we are ready to apply our main theorem \ref{main}! We know from previous subsection that there's a function $f(a,b,h)$ of at most degree $4$ such that a triangle satisfies Euler's theorem if and only if $f(a,b,h)=0$, moreover we have just discovered that if $(a,b,h)$ satisfy one of the equations \eqref{iso}
then also $f(a,b,h)=0$ since the theorem is true. Taking
$$P(a,b,h)=(a^2-b^2-h^2)(a-2b)(h^2+b^2-2ab)$$
we only need to find a line which has $5$ single intersections with $P(a,b,c)$, in order to apply \ref{main} to do the job, proving that $f=0$ everywhere.
Fortunately, this line $\vec s: \R\to \R^3$ is not hard to find, take for example 

$$\vec s(k)=\begin{pmatrix} 2 \\ k \\ 1 \end{pmatrix};\ \ \ \ \ \ \ \ \ \begin{pmatrix} a=2 \\ b=k \\ h=1 \end{pmatrix}$$
which gives some simple intersections:
$$(a^2-b^2-h^2)=0\implies 3-k^2=0$$
$$a-2b=0 \implies k=1$$
$$(h^2+b^2-2ab)=0 \implies k^2-4k+1=0$$
which are all distinct: $k=\{-\sqrt 3,\ \sqrt 3,\ 1,\ 2-\sqrt 3,\ 2+\sqrt 3\}$. This implies with \ref{main} that $f(a,b,h)=0$ for any $a,b,h$ i.e. the Euler line exists for every triangle.

%%%%%%%%%%%%%%%%%%%%%%%%%%%%%%%%%%%%%%%%%%%%%%%
% ZONA appestati, non passare
%%%%%%%%%%%%%%%%%%%%%%%%%%%%%%%%%%%%%%%%%%%%%%%

\appendix{
\section{Simple results to estimate the degree of the associated polynomial}\label{boring}

\begin{lem}
Let $A=\vec f_A(a,b,c)$ and $B=\vec f_B(a,b,c)$ two points whose coordinates are polynomial functions of degree $g$ of $a,b,c$. Then, the line passing from $A$ and $B$ has the form\footnote{$\vec r$ is a function of $k$ which implicitly depends on $a,b,c$ as all the geometric elemnts we are considering, next this implicit dependence will be omitted}
$$\vec r(k;\ a,b,c)=\vec v_r(a,b,c)k+\vec q_r (a,b,c)$$
where $\vec v_r(a,b,c)$ and $\vec q_r (a,b,c):\ \R^3\to \R^2$ have degree at most degree $g$.
\label{2ptR}
\end{lem}
\begin{proof}
We can just choose $\vec v_r=\vec f_B(a,b,c)-\vec f_A(a,b,c)$ and $\vec q_r=\vec f_A(a,b,c)$ and verify that
$$\vec r(0)=\vec f_A(a,b,c)=A$$
$$\vec r(1)=\vec f_B(a,b,c)=B$$
then, $\vec v_r,\vec q_r$ have at most degree $p$ being linear combinations of $\vec f_A,\vec f_B$.
\end{proof}

\begin{lem}
Let 
$$\vec r(k;\ a,b,c)=\vec v_r(a,b,c)k+\vec q_r (a,b,c)$$
be a line in $\R^2$ such that $\vec v_r,\vec q_r: \R^3\to \R^2$ are linear functions of $a,b,c$ and $\ell(a,b,c): \R^3\to \R$ be another linear function.
Then intersection point $P$ of $\vec r(k)$ and the vertical axis $\{(x,y): x=\ell(a,b,c)\}$
\label{2Rpt} is given by the function $\vec f_P(a,b,c): \R^3 \to \R^2$ whose components are 
$$f_{P,x}(a,b,c)=\ell(a,b,c)$$
$$f_{P,y}(a,b,c)$$
where $f_{P,y}(a,b,c)$ is a rational function of $a,b,c$, where the numerator has at most degree 2 and the denominator degree 1.
\end{lem}

\begin{proof}
To find the intersection, it is enough to impose that $r_x(k)$ is equal to $\ell(a,b,c)$:
$$k=\frac{\ell(a,b,c)-q_{r,x}(a,b,c)}{v_{r,x}(a,b,c)}$$
where, by assumpion, the numerator and the denominator are first degree on $(a,b,c)$.
Now we can substitute in the equation of $\vec r(k)$ to find the intersection point
$$P:=\begin{pmatrix} f_{P,x}(a,b,c)\\
f_{P,y}(a,b,c)\end{pmatrix}$$
$$f_{P,x}(a,b,c)=\ell(a,b,c)$$

$$f_{P,y}(a,b,c)=\frac{\ell(a,b,c)-q_{r,x}(a,b,c)}{v_{r,x}(a,b,c)}    \vec v_r(a,b,c)+\vec q_r(a,b,c)$$
where the second component has, by construction, has at most degree $2$(numerator) on $1$(denominator)
\end{proof}

you can imagine as many as you wish of lemmas line this! For example, you could prove that
\begin{exe}
Let
$$\vec r(k;\ a,b,c)=\vec v_r(a,b,c)k+\vec q_r (a,b,c)$$
$$\vec s(k;\ a,b,c)=\vec v_s(a,b,c)k+\vec q_s (a,b,c)$$
be two lines in $\R^2$ such that $\vec r,\vec s$ are linear in $(a,b,c)$. Then their intersection point $P=\vec f_P(a,b,c)$, if exists, has components
which are rational functions of $a,b,c$ where the numerator has at most degree 3 and the denominator degree 2.
\end{exe}
The argument is exactly the same and it reduces just to a simple exercise of basic algebra.
}
ok

\end{document}